\documentclass[12pt,leqeq]{article}

\usepackage{amssymb}

\usepackage{amsmath}

\usepackage[all]{xy}

 \newcommand{\Z}{\mathbb{Z}}
 \newcommand{\Agr}{A^\bullet}
 
 \newcommand{\Mgr}{M^\bullet}

 \newcommand{\kr}{\mathbf{k}}
 
 \newcommand{\N}{\mathbb{N}}
 \newcommand{\Real}{\mathbb{R}}
 \newcommand{\Complex}{\mathbb{C}}

\newtheorem{thm}{Theorem}

\newtheorem{defn}[thm]{Definition}

\begin{document}

\setlength{\baselineskip}{16pt}

\title{On the $q$-analogue of the Maurer-Cartan equation}

\author{Mauricio Angel\hspace{.3cm}and Rafael D\'\i az}

\maketitle

\begin{sloppypar}

\begin{abstract}
We consider deformations of the differential of a $q$-differential
graded algebra. We prove that it is controlled by a generalized
Maurer-Cartan equation. We find explicit formulae for the
coefficients $c_k$ involved in that equation.
\end{abstract}

\section*{Introduction}

Following the work of Kapranov \cite{Kap} the interest on the
theory of N-complexes has risen. A delicate issue in the theory is
the appropriate definition of a $N$-differential algebra. At the
moment there exists two options for such a definition, each
describing very distinct mathematical phenomena.

\smallskip

In \cite{AD} the notion of $N$-differential graded algebra $A$ has
been defined as follows: $A$ must be a graded associative algebra
provided with an operator $d:A\to A$ of degree 1 such that
$d(ab)=d(a)b+(-1)^{\bar{a}}ad(b)$ and $d^N=0$.

\smallskip

The study of the deformations of the differential of a N-dga has
also been initiated in \cite{AD}. Applications to differential
geometry are given in \cite{AD2}. Applications to the study of
$A_\infty$-algebras of depth $N$, are contained in \cite{AD3}.

\smallskip

The other possible $N$-generalization of the notion of a
differential algebra (see \cite{KaWa},\cite{KN},\cite{Sit}) goes
as follows: One fixed a primitive $N$-th root of unity $q$. One
considers $q$-differential graded algebras $A$ defined by: $A$
must be a graded associative algebra provided with an operator
$d:A\to A$ of degree 1 such that $d(ab)=d(a)b+q^{\bar{a}}ad(b)$
and $d^N=0$.

\smallskip

In this paper we consider deformations of the differential of a
$q$-differential graded algebra. We prove that it is controlled by
a generalized Maurer-Cartan equation. We find explicit formulae
for the coefficients $c_k$ involved in that equation.

\section{$q$-differential graded algebras and modules}
Let $\kr$ be a commutative ring with unity. For $q\in \kr$ the
$q$-numbers are given by $[\ ]_q:\N\to\kr$, $k\longmapsto
[k]_q=1+q+\cdots q^{k-1}$ and $[k]_q!=[1]_q[2]_q\cdots[k]_q$. We
assume that on $\kr$, for $N\in \N$ fixed, we can choose $q\in\kr$
such that $[N]_q=0$ and $[k]_q$ is invertible for $1\leq k\leq
N-1$, in this case we refer to $q$ as a primitive $N$-th root of
unity.

\begin{defn}\label{qdga}
A {\bf $q$-differential graded algebra} or q-dga over $\kr$ is a
triple $(\Agr,m,d)$ where $m:A^k\otimes A^l\to A^{k+l}$ and
$d:A^k\to A^{k+1}$ are $\kr$-modules morphisms satisfying
\begin{enumerate}
\item[1)] $(\Agr,m)$ is a graded associative $\kr$-algebra.

\item[2)] $d$ satisfies
the q-Leibniz rule $d(ab)=d(a)b+q^{\bar{a}}ad(b)$.

\item[3)] If $q$ is a primitive $N$-th root of unity, then
$d^{N}=0$, i.e., $(\Agr,d)$ is a $N$-complex.

\end{enumerate}

\end{defn}

\begin{defn}
Let $(\Agr,m_A,d_A)$ be a q-dga and $\Mgr$ be a $\kr$-module. A
{\bf $q$-differential module} (q-dgm) over $(\Agr,m_A,d_A)$ is a
triple $(\Mgr,m_M,d_M)$ where $m_M:A^k\otimes M^l\to M^{k+l}$ and
$d_M:M^k\to M^{k+1}$ are $\kr$-module morphisms satisfying the
following properties
\begin{enumerate}
\item $m_M(a,m_M(b,\phi))=m_M(m_A(a,b),\phi)$, for all
$a,b\in\Agr$ and $\phi\in\Mgr$. We denote $m_M(a,\phi)$ by
$a\phi$, for all $a\in A$ and $\phi\in M$.

\item $d_M(a\phi)=d_A(a)\phi+q^{\bar{a}}ad_M(\phi)$ for all
$a\in\Agr$ and $\phi\in\Mgr$.

\item If $q$ is a primitive $N$-th root of unity, then
$d_M^{N}=0$.
\end{enumerate}
\end{defn}

\section{N-curvature}

Let $(\Mgr,m_M,d_M)$ be a $q$-dgm over a $q$-dga $(\Agr,m_A,d_A)$.
For $a\in M^1$ we define a deformation of $d_M$ as follows
\[D\phi=d_M\phi+a\phi,\hspace{.5cm}\text{for all $\phi\in\Mgr$}.\]
In order to obtain a $q$-dgm structure on $(\Mgr,m_M,D)$ we
require that $D^N=0$. In [KN] the consecutive powers of $D$ are
calculated and the following formula is obtained
\[D^{N}\phi=d_M^{N}\phi+\sum_{k=1}^{N-1}\binom{N}{k}_q(D^{k-1}a)d^{N-k}\phi+(D^{N-1}a)\phi,\]
where $\binom{N}{k}_q$ is the $q$-binomial coefficients given by
\[\binom{N}{k}_q=\frac{[N]_q!}{[N-k]_q![k]_q!}.\]
For $q$ a primitive N-th root of unity the formula above reduces
to
\begin{equation}\label{eq1}
D^N\phi=(D^{N-1}a)\phi,
\end{equation}
because in this case $d_M^N=0$ and $\binom{N}{k}_q=0$ for $0<k<N$,
since $[N]_q=1+q+\cdots+q^{N-1}=\frac{1-q^N}{1-q}=0$.

\smallskip

Formula $(\ref{eq1})$ is the first step towards the construction
of the $q$-analogue of the Maurer-Cartan equation. A further step
is required in order to write down $D^{N-1}a$ in terms of $d_M$
and $a$ only. We shall work with a slightly more general problem:
finding an explicit expression for $D^N\phi$ given in terms of
$d_M$ and $a$ only. First we review some notations from [AD].

\smallskip

For $s=(s_1,...,s_n)\in \N^n$ we set $l(s)=n$ and
$|s|=\sum_i{s_i}$. For $1\leq i <n,\ s_{>i}$ denotes the vector
given by $s_{>i}=(s_{i+1},...,s_n)$, for $1<i\leq n,\ s_{<i}$
stands for $s_{<i}=(s_1,...,s_{i-1})$, we also set
$s_{>n}=s_{<1}=\emptyset$. $\N^{(\infty)}$ denotes the set
$\bigsqcup_{n=0}^{\infty}\N^n$. By convention
$\N^{0}=\{\emptyset\}$.

\smallskip

For $e\in End(\Mgr)$ and $s\in\N^n$ we define
$e^{(s)}=e^{(s_1)}...e^{(s_n)}$, where $e^{(l)}=d_{End}^{l}(e)$ if
$l\geq 1$, $e^{(0)}=e$ and $e^{\emptyset}=1$. In the case that
$e_a\in End(\Mgr)$ is given by
\[e_a(\phi)=a\phi,\hspace{.3cm}\text{for $a\in M^1$ fixed and all $\phi\in\Mgr$,}\]
then $e_a^{(l)}=d_{End}^{l}(e_a)$ reduces to
$e_a^{(l)}=e_{d^{l}(a)}$, thus
\[e_a^{(s)}=e_a^{(s_1)}\cdots e_a^{(s_n)}=e_{d^{s_1}(a)}\cdots e_{d^{s_n}(a)}.\]

For $N\in\N$ we define $E_N=\{s\in\N^{(\infty)}:|s|+l(s)\leq N\}$
and for $s\in E_N$ we define $N(s)\in\Z$ by $N(s)=N-|s|-l(s)$.

\smallskip

We introduce a discrete quantum mechanical \footnote{a discrete
quantum mechanical system is given by the following data (1) A
directed graph $\Gamma$ (finite or infinite). (2) A map
$L:E_\Gamma\to\Real$ called the Lagrangian map of the system. An
associated Hilbert space $\mathcal{H}=\Complex^{V_\Gamma}$.
Operators $U_n:\mathcal{H}\to\mathcal{H}$, where $n\in\Z$, given
by
\[(U_nf)(y)=\sum_{x\in V_\Gamma}\omega_n(y,x)f(x),\]
where the discretized kernel $\omega_n(y,x)$ admits the following
representation
\[\omega_n(y,x)=\sum_{\gamma\in P_n(\Gamma,x,y)}\prod_{e\in\gamma}e^{iL(e)}.\]
$P_n(\Gamma,x,y)$ denotes the set of length $n$ paths in $\Gamma$
from $x$ to $y$, i.e., sequences $(e_1,\cdots,e_n)$ of edges in
$\Gamma$ such that $s(e_1)=x$, $t(e_i)=s(e_{i+1}),\ i=1,...,n-1$
and $t(e_n)=y$.} by
\begin{enumerate}
\item $V=\N^{(\infty)}$. \item There is a unique directed edge
$e$ from vertex $s$ to $t$ if and only if
$t\in\{(0,s),s,(s+e_i)\}$ where
$e_i=(0,..,\underset{\scriptsize{i-th}}{\underbrace{1}},..,0)\in\N^{l(s)}$,
in this case we set $source(e)=s$ and $target(e)=t$.

\item Edges are weighted according to the following table
\smallskip
\begin{center}
\begin{tabular}{|l|l|l|}\hline $source(e)$    &      $target(e)$ &
$v(e)$      \\  \hline
$s$         &      $(0,s)$      & $1$                 \\
$s$         &      $s$          & $q^{|s|+l(s)}$    \\
$s$         &     $(s+e_i)$    & $q^{|s_{<i}|+i-1}$ \\
\hline
\end{tabular}
\end{center}
\end{enumerate}

The set $P_N(\emptyset,s)$ consists of all paths
$\gamma=(e_1,...,e_N)$, such that $source(e_1)=\emptyset$,
$target(e_N)=s$ and $source(e_{l+1})=target(e_l)$. For $\gamma\in
P_N(\emptyset,s)$ we define the weight $v(\gamma)$ of $\gamma$ as
\[v(\gamma)=\prod_{l=1}^{N}v(e_l).\]

The following result is proven as Theorem 17 in [AD].
\begin{thm}\label{MCMN}
Let $(\Mgr,m_M,d_M)$ be a q-dgm over a q-dga $(\Agr,m_A,d_A)$. For
$a\in M^1$ consider the map $D\phi=d_M\phi+a\phi$. Then we have
\[D^N=\sum_{k=0}^{N-1}c_k d_M^k\]
where
\[c_k=\sum_{\begin{subarray}{c} s\in E_N\\ N(s)=k \\ s_i<N\\ \end{subarray}}c_q(s,N)a^{(s)}
\hspace{.5cm}\text{and}\hspace{.5cm} c_q(s,N)=\sum_{\gamma\in
P_N(\emptyset,s)}v(\gamma).\]
\end{thm}

\smallskip

Suppose that we have a $q$-dgm $(\Mgr,m_M,d_M)$, $q$ a 3-rd
primitive root of unity, and we want to deform it to $D=d_M+a$,
$a\in M^1$. So we required that $D^3=0$. By Theorem \ref{MCMN} we
must have $\displaystyle{\sum_{k=0}^{2}c_k d_M^k=0}$. Let us
compute the coefficients $c_k$. First, notice that
\[E_3=\{\emptyset, (0), (1), (2),(0,0), (1,0), (0,1), (0,0,0)\}\]
We proceed to compute the coefficients $c_k$ for $0\leq k\leq 2$
using Theorem \ref{MCMN}.

\smallskip

$\underline{k=0}$
\smallskip

There are four vectors in $E_3$ such that $N(s)=0$, these are
$(2), (1,0), (0,1), (0,0,0)$.

\smallskip

For $s=(2)$ the only path from $\emptyset$ to $(2)$ of length 3 is
$\emptyset\to(0)\to(1)\to(2)$\hspace{.3cm} with weight 1 and
$e_a^{(2)}=d_M^2(e_a)$, thus $c(s,3)=d_M^2(a)$.

\medskip

For $s=(1,0)$ the only path from $\emptyset$ to $(1,0)$ of length
3 is $\emptyset\to(0)\to(0,0)\to(1,0)$\hspace{.3cm} with weight 1
and $e_a^{(1,0)}=d_M(e_a)e_a$, thus $c(s,3)=d_M(a)e_a$.

\medskip

For $s=(0,1)$ the paths from $\emptyset$ to $(0,1)$ of length 3
are $\emptyset\to(0)\to(0,0)\to(0,1)$\hspace{.3cm}with weight $q$;
and $\emptyset\to(0)\to(1)\to(0,1)$\hspace{.3cm}with weight 1.
Since $e_a^{(0,1)}=e_ad_M(e_a)$, then $c(s,3)=(1+q)e_ad_M(e_a)$.

\medskip

For $s=(0,0,0)$ the only path from $\emptyset$ to $(0,0,0)$ of
length 3 is $\emptyset\to(0)\to(0,0)\to(0,0,0)$\hspace{.3cm} whose
weight is 1 and $e_a^{(0,0,0)}=a^3$, thus $c(s,3)=a^3$.

\bigskip

$\underline{k=1}$

\smallskip

There are two vectors in $E_3$ such that $N(s)=1$, namely $(1)$
and $(0,0)$.

\smallskip

For $s=(1)$, the paths from $\emptyset$ to $(1)$ of length 3 are
$\emptyset\to\emptyset\to(0)\to(1)$\hspace{.3cm}with weight 1;
$\emptyset\to(0)\to(0)\to(1)$\hspace{.3cm}with weight $q$;
$\emptyset\to(0)\to(1)\to(1)$\hspace{.3cm}with weight $q^2$.
$e_a^{(1)}=d_M(e_a)$ and thus $c_q(s,3)=(1+q+q^2)=0$.

\medskip

For $s=(0,0)$ the paths from $\emptyset$ to $(0,0)$ of length 3
are $\emptyset\to(0)\to(0)\to(0,0)$\hspace{.3cm}with weight $q$.
$\emptyset\to\emptyset\to(0)\to(0,0)$\hspace{.3cm}with weight 1.
$\emptyset\to(0)\to(0,0)\to(0,0)$\hspace{.3cm}with weight $q^2$.
$e_a^{(0,0)}=a^2$ and $c_q(s,3)=0$.

\bigskip

$\underline{k=2}$

\smallskip

$(0)$ is the only vector in $E_3$ such that $N(s)=2$.

\smallskip

The paths from $\emptyset$ to $(0)$ of length 3 are
$\emptyset\to\emptyset\to\emptyset\to(0)$\hspace{.3cm}with weight
1; $\emptyset\to\emptyset\to(0)\to(0)$\hspace{.3cm}with weight
$q$; $\emptyset\to(0)\to(0)\to(0)$\hspace{.3cm}with weight $q^2$.
$e_a^{(0)}=a$ and $c(s,3)=(1+q+q^2)=0$.

\medskip

So we have proven that the 3-curvature is given by
\[D^3=d_M^2(a)+d_M(a)a+(1+q)ad_M(a).\]
Notice that in case of $q$-commutativity
$\alpha\beta=q^{\bar{\alpha}\bar{\beta}}\beta\alpha$, the
3-curvature reduces to $D^3=d_M^2(a)$.

\bigskip

Suppose now that $q$ is a 4-th primitive root of unity and we
required that $D^4=0$. By Theorem \ref{MCMN} we must have
$\displaystyle{\sum_{k=0}^{3}c_k d_M^k=0}$. Notice that
\[E_4=\{\emptyset, (0), (1), (2), (3),(0,0), (1,0), (0,1), (2,0), (0,2), (1,1), \]
\[(0,0,0), (1,0,0), (0,1,0), (0,0,1),(0,0,0,0)\}\]
We proceed to compute the coefficients $c_k$ for $0\leq k\leq 3$
using Theorem \ref{MCMN}.

\smallskip

$\underline{k=3}$\\
$\emptyset\stackrel{i}{\longrightarrow}\emptyset\stackrel{\longrightarrow}{}
(0)\stackrel{j}{\longrightarrow}(0)$\hspace{.3cm} the weight is
$\displaystyle{\sum_{i+j=3}q^{j}}$ and $c_3=(1+q+q^2+q^3)a=0$.

$\underline{k=2}$\\
$\emptyset\stackrel{i}{\longrightarrow}\emptyset\longrightarrow(0)\stackrel{j}{\longrightarrow}(0)\longrightarrow(0,0)\stackrel{k}{\longrightarrow}(0,0)$\hspace{.2cm}
the weight is $\displaystyle{\sum_{i+j+k=2}q^{j}q^{2k}}$ and
$c_q((0,0),4)a^{(0,0)}=(1+q^2)a^2$.

\smallskip

\hspace{.4cm}$\emptyset\stackrel{i}{\longrightarrow}\emptyset\longrightarrow(0)\stackrel{j}{\longrightarrow}(0)\longrightarrow(1)\stackrel{k}{\longrightarrow}(1)$\hspace{.2cm}
the weight is $\displaystyle{\sum_{i+j+k=2}q^{j}q^{2k}}$ and
$c_q((1),4)a^{(1)}=(1+q^2)d(a)$.

\medskip

Finally, $c_2=(1+q^2)(a^2+d(a))$.

$\underline{k=1}$\\
$\emptyset\stackrel{i}{\longrightarrow}\emptyset\longrightarrow(0)\stackrel{j}{\longrightarrow}(0)\longrightarrow(0,0)\stackrel{k}{\longrightarrow}(0,0)\longrightarrow(0,0,0)\stackrel{l}{\longrightarrow}(0,0,0)$\\
the weight is $\displaystyle{\sum_{i+j+k+l=1}q^{j}q^{2k}q^{3l}}$
and $c_q((0,0,0),4)a^{(0,0,0)}=(1+q+q^2+q^3)a^3=0$.
\[\emptyset\stackrel{i}{\longrightarrow}\emptyset\longrightarrow(0)\stackrel{j}{\longrightarrow}(0)\longrightarrow(1)\stackrel{k}{\longrightarrow}(1)\longrightarrow(0,1)\stackrel{l}{\longrightarrow}(0,1)\]
the weight is $\displaystyle{\sum_{i+j+k+l=1}q^{j}q^{2k}q^{3l}}$.
\[\emptyset\stackrel{i}{\longrightarrow}\emptyset\longrightarrow(0)\stackrel{j}{\longrightarrow}(0)\longrightarrow(0,0)\stackrel{k}{\longrightarrow}(0,0)\longrightarrow(0,1)\stackrel{l}{\longrightarrow}(0,1)\]
the weight is $\displaystyle{\sum_{i+j+k+l=1}q^{j}q^{2k}qq^{3l}}$
and $c_q((0,1),4)a^{(0,1)}=((1+q+q^2+q^3)+(1+q+q^2+q^3))ad(a)=0$.

\[\emptyset\stackrel{i}{\longrightarrow}\emptyset\longrightarrow(0)\stackrel{j}{\longrightarrow}(0)\longrightarrow(0,0)\stackrel{k}{\longrightarrow}(0,0)\longrightarrow(1,0)\stackrel{l}{\longrightarrow}(1,0)\]
the weight is $\displaystyle{\sum_{i+j+k+l=1}q^{j}q^{2k}q^{3l}}$
and $c_q((1,0),4)a^{(1,0)}=(1+q+q^2+q^3)d(a)a=0$.

\[\emptyset\stackrel{i}{\longrightarrow}\emptyset\longrightarrow(0)\stackrel{j}{\longrightarrow}(0)\longrightarrow(1)\stackrel{k}{\longrightarrow}(1)\longrightarrow(2)\stackrel{l}{\longrightarrow}(2)\]
the weight is $\displaystyle{\sum_{i+j+k+l=1}q^{j}q^{2k}q^{3l}}$
and $c_q((2),4)a^{(2)}=(1+q+q^2+q^3)d^2(a)=0$.

Finally,
$c_1=c_q((0,0,0),4)a^{(0,0,0)}+c_q((0,1),4)a^{(0,1)}+c_q((1,0),4)a^{(1,0)}+c_q((2),4)a^{(2)}=0$.

$\underline{k=0}$\\
$\emptyset\longrightarrow(0)\longrightarrow(0,0)\longrightarrow(0,0,0)\longrightarrow(0,0,0,0)$\\
the weight is $1$ and $c_q((0,0,0,0),4)a^{(0,0,0,0)}=a^4$.

\[\emptyset\longrightarrow(0)\longrightarrow(0,0)\longrightarrow(0,0,0)\longrightarrow(0,0,1)\]
\[\emptyset\longrightarrow(0)\longrightarrow(1)\longrightarrow(0,1)\longrightarrow(0,0,1)\]
\[\emptyset\longrightarrow(0)\longrightarrow(0,0)\longrightarrow(0,1)\longrightarrow(0,0,1)\]
the weight is $1+q+q^2$ and
$c_q((0,0,1),4)a^{(0,0,1)}=(1+q+q^2)a^2d(a)$.

\[\emptyset\longrightarrow(0)\longrightarrow(0,0)\longrightarrow(1,0)\longrightarrow(0,1,0)\]
\[\emptyset\longrightarrow(0)\longrightarrow(0,0)\longrightarrow(0,0,0)\longrightarrow(0,1,0)\]
the weight is $1+q$ and $c_q((0,1,0),4)a^{(0,1,0)}=(1+q)ad(a)a$.

\[\emptyset\longrightarrow(0)\longrightarrow(0,0)\longrightarrow(0,0,0)\longrightarrow(1,0,0)\]
the weight is $1$ and $c_q((1,0,0),4)a^{(1,0,0)}=d(a)a^2$.

\[\emptyset\longrightarrow(0)\longrightarrow(0,0)\longrightarrow(1,0)\longrightarrow(2,0)\]
the weight is $1$ and $c_q((2,0),4)a^{(2,0)}=d^2(a)a$.

\[\emptyset\longrightarrow(0)\longrightarrow(0,0)\longrightarrow(0,1)\longrightarrow(0,2)\]
\[\emptyset\longrightarrow(0)\longrightarrow(1)\longrightarrow(0,1)\longrightarrow(0,2)\]
\[\emptyset\longrightarrow(0)\longrightarrow(1)\longrightarrow(2)\longrightarrow(0,2)\]
the weight is $1+q+q^2$ and
$c_q((0,2),4)a^{(0,2)}=(1+q+q^2)ad^2(a)$.

\[\emptyset\longrightarrow(0)\longrightarrow(0,0)\longrightarrow(1,0)\longrightarrow(1,1)\]
\[\emptyset\longrightarrow(0)\longrightarrow(1)\longrightarrow(0,1)\longrightarrow(1,1)\]
the weight is $1+q^2$ and $c_q((1,1),4)a^{(1,1)}=(d(a))^2$.

\[\emptyset\longrightarrow(0)\longrightarrow(1)\longrightarrow(2)\longrightarrow(3)\]
the weight is $1$ and $c_q((3),4)a^{(3)}=d^3(a)$.

Finally,
$c_0=a^4+(1+q+q^2)a^2d(a)+(1+q)ad(a)a+d(a)a^2+d^2(a)a+(1+q+q^2)ad^2(a)+(d(a))^2+d^3(a)$.

\medskip

So we have that the 4-curvature is given by
\[D^4=(1+q^2)(a^2+d(a))d_M^2+ a^4+(1+q+q^2)a^2d(a)+(1+q)ad(a)a+\]
\[d(a)a^2+d^2(a)a+(1+q+q^2)ad^2(a)+(d(a))^2+d^3(a).\]

\bigskip

In the remainder of the paper we consider infinitesimal
deformations of a $q$-differential.

\begin{thm}
Let $(\Mgr,m_M,d_M)$ be a $q$-dgm such that $d_M^N=0$, consider
the infinitesimal deformation $D=d_M+te$, where $e\in End^1(\Mgr)$
and $t^2=0$, then
\[D^N=\sum_{k=0}^{N-1}\left(\sum_{v\in Par(N,N-k-1)}q^{|v|+w(v)}\right)d_M^{N-k-1},\]
where
\[Par(N,N-k-1)=\{(v_0,\cdots,v_{N-k-1})/ \sum v_i=N\},\]
$|v|=v_1+\cdots +v_n$ and $w(k)=\sum iv_i$.
\end{thm}
{\bf Proof}: By Theorem \ref{MCMN},
$D^N=\sum_{k=0}^{N-1}c_kd_M^k$. Since $t^2=0$, then
\[(te)^{(s)}=(te)^{(s_1)}\cdots(te)^{(s_{l(s)})}=t^{l(s)}e^{(s)}=0\hspace{.3cm}\text{unless $l(s)\leq 1$}.\]
Thus $E_N$ is given by
\[E_N=\{(0),(1),\cdots,(N-1)\}.\]
For $0\leq k\leq N-1$, since $N(s)=N-|s|-l(s)=k$ and $l(s)=1$ we
obtain $|s|=N-k-1$, and the only vector $s$ in $E_N$ of length 1
such that $|s|=N-k-1$ is precisely $s=(N-k-1)$. Thus
\[c_k=\sum_{\begin{subarray}{c} s\in E_N\\ N(s)=k \\ s_i<N\\
\end{subarray}}c_q(s,N)a^{(s)}=c_q((N-k-1)),N)a^{(s)}=c_q((N-k-1)),N)d^{N-k-1}(a).\]

\smallskip

Any path from $\emptyset$ to $(N-k-1)$ of length $N$ must be of
the form
\[\underbrace{\emptyset\rightarrow\cdots\rightarrow\emptyset}_{p_0}\rightarrow\underbrace{(0)\rightarrow\cdots\rightarrow(0)}_{p_1}\rightarrow\underbrace{(1)\rightarrow\cdots(1)}_{p_2}\rightarrow\cdots\underbrace{(N-k-1)\rightarrow\cdots\rightarrow(N-k-1)}_{p_{N-k-1}}\]
with $p_0+p_1+\cdots+p_{N-k-1}=N$. The weight of such path is
$q^{p_1}q^{p_2(1+1)}...q^{p_{N-k-1}(N-k)}=q^{|p|+w(p)}$, where
$p=(p_0,\cdots,p_{N-k-1})$. Then
\[c_q(s,N)=\sum_{\gamma\in P_N(\emptyset,s)}v(\gamma)=\sum_{p\in Par(N,N-k-1)}q^{|p|+w(p)}.\blacklozenge\]


\subsection*{Acknowledgement}
We thank Sylvie Paycha and Nicol\'as Andruskiewitsch.


\[\begin{array}{c}

\mbox{Mauricio Angel. Universidad Central de Venezuela (UCV).} \ \  \mbox{\texttt{mangel@euler.ciens.ucv.ve}} \\

\mbox{Rafael D\'\i az. Universidad Central de Venezuela (UCV).} \ \  \mbox{\texttt{rdiaz@euler.ciens.ucv.ve}} \\
\end{array}\]

\end{sloppypar}
\end{document}